\documentclass[12pt]{article}
\usepackage{graphicx}
\usepackage{amsmath, amssymb,amsfonts}

\newcommand{\dint}{\displaystyle \int}
\renewcommand{\dfrac}{\displaystyle \frac}

\newtheorem{theorem}{Theorem}[section]
\newtheorem{lemma}[theorem]{Lemma}
\newtheorem{proposition}[theorem]{Proposition}
\newtheorem{remark}[theorem]{Remark}

\newtheorem{definition}[theorem]{Definition}

\begin{document}

\title{\textbf{Stability and genericity for SPDEs \\ driven by spatially correlated noise}}
\date{ \ \ }
\author{}
\maketitle

\begin{center}
\vspace{-2.5cm}

{\small \sc Khaled Bahlali} \\
{\small D\'{e}partement de Math\'{e}matiques, Universit\'{e} du Sud, Toulon--Var} \\
{\small B.P. 132, 83957 La Garde, France, bahlali@univ-tln.fr}\\

\medskip 

{\small \sc M'hamed Eddahbi\footnote[1]{Corresponding author}}\\
{\small D\'{e}partement de Math. \& Info. FSTG, Universit\'{e} Cadi Ayyad}\\
{\small B.P. 549, Marrakech, Maroc, eddahbi@fstg-marrakech.ac.ma} \\

\medskip 

{\small \sc Mohamed Mellouk}\\
{\small Universit\'{e} Montpellier 2, Institut de Math\'{e}matiques et de Mod\'{e}lisation}\\
{\small de Montpellier, F-34095 Montpellier Cedex 5, France, mellouk@math.univ-montp2.fr}
\end{center}
\vspace{0,2cm}
\begin{abstract}
We consider stochastic partial differential equations  on
$\mathbb{R}^{d},\, d\geq 1$, driven by a Gaussian noise white in
time and colored in space, for which the pathwise uniqueness holds.
By using the Skorokhod representation theorem we establish various
strong stability results. Then, we give an application to the
convergence of the Picard successive approximation. Finally, we show that in
the sense of Baire category, almost all stochastic partial
differential equations with continuous and bounded coefficients have
the properties of existence and uniqueness of solutions as well as
the continuous dependence on the coefficients.

\vspace{2cm} 


\medskip


\medskip

\medskip

\noindent \textbf{Acknowledgements.} This work is partially supported by l'Action Int\'{e}gr\'{e}e  MA/142/06.

\end{abstract}

\newpage 

\section{Introduction and general framework}

The paper in concerned with stochastic partial differential equations
(SPDEs) of the form

\begin{eqnarray}\label{eq0}
\left\{ 
\begin{array}{l}
Lu(t,x)=\sigma (t,x,u(t,x)){\dot{F}}(t,x)+b(t,x,u(t,x)), \\ \\
u(0,x)=0,\;\dfrac{\partial u}{\partial t}(0,x)=0,
\end{array}
\right. 
\end{eqnarray}
where, $t\in [0,T]$ for some fixed $T>0$, $x\in \mathbb{R}^{d},d\geq 1$
and $L$ is a second order partial differential operator. The coefficients $%
\sigma $ and $b:\mathbb{R}\rightarrow \mathbb{R}$ are given measurable functions.
Our SPDEs include, for instance, the stochastic heat and wave equations in
spatial dimension $d\geq 1$. For the simplicity, we will assume that the
initial condition is null. The result can, however, be properly extended to
cover smooth initial conditions.\newline

\noindent Let $\mathcal{D}(\mathbb{R}^{d+1})$ be the space of all infinitely
differentiable fonctions with compact support. On a probability space $%
(\Omega ,\mathcal{G},P)$, the noise $F=\{F(\varphi ),\varphi \in \mathcal{D}(%
\mathbb{R}^{d+1})\}$ is assumed to be an $L^{2}(\Omega ,\mathcal{G},P)$--valued
Gaussian process with mean zero and covariance functional given by 
\[
J(\varphi ,\psi )=\int_{\mathbb{R}_{+}}ds\int_{\mathbb{R}^{d}}\Gamma (dx)\left(
\varphi (s,\cdot )\ast \widetilde{\psi }(s,\cdot )\right) (x),
\]
where $\widetilde{\psi }(s,x)=\psi (s,-x)$ and $\Gamma $ is a non--negative
and non--negative definite tempered measure, therefore symmetric. Let $\mu $
denote the spectral measure of $\Gamma $, which is also a tempered measure.%
\newline

\noindent Denote by $\mathcal{F}\varphi $ the Fourier transform of $\varphi $%
. Clearly $\mu =\mathcal{F}^{-1}(\Gamma )$. This gives 
\[
J(\varphi ,\psi )=\int_{\mathbb{R}_{+}}ds\int_{\mathbb{R}^{d}}\mu (d\xi )\mathcal{F%
}\varphi (s,\cdot )(\xi )\overline{\mathcal{F}\psi (s,\cdot )}(\xi ),
\]
where $\overline{z}$ is the complex conjugate of $z.$\newline

\noindent Following the same approach in \cite{Da}, the Gaussian process $F$
can be extended to a worthy martingale measure $M=\{M(t,A):= F([0,t]\times
A)\,:\,t\geq 0, \, A\in \mathcal{B}_{b}(\mathbb{R}^{d})\}$ which shall acts as
integrator, in the Walsh sence \cite{W}, where $\mathcal{B}_{b}(\mathbb{R}^{d})$
denotes the bounded Borel subsets of $\mathbb{R}^d$. Let $\mathcal{G}_{t}$ be
the completion of the $\sigma $--field generated by the random variables $%
\{M(s,A),\;0\leq s\leq t,\;A\in \mathcal{B}_{b}(\mathbb{R}^{d})\}$. The
properties of $F$ ensure that the process $M=\{M(t,A),\;t\geq 0,\;A\in 
\mathcal{B}_{b}(\mathbb{R}^{d})\}$, is a martingale with respect to the
filtration $\{\mathcal{G}_{t}:t\geq 0\}$. \newline

\noindent One can give a rigorous meaning to solution of equation (\ref{eq0}%
), by means of a jointly measurable and $\mathcal{G}_{t}$--adapted process $%
\{u(t,x):(t,x)\in \mathbb{R}_{+}\times \mathbb{R}^{d}\}$ satisfying, for each $%
t\geq 0$ and $x\in \mathbb{R}^{d}$, a.s. the following evolution equation : 
\begin{eqnarray}  \label{integrale}
u(t,x) &=&\int_{0}^{t}\int_{\mathbb{R}^{d}}S(t-s,x,y)\sigma (s,y,u(s,y))M(ds,dy)
\nonumber  \label{eq1} \\
&&\,\,\,\,\,\,\,+\,\int_{0}^{t}ds\int_{\mathbb{R}^{d}}dyS(t-s,x,y)b(s,y,u(s,y))
\end{eqnarray}
where $S(t,x,y)$ stands for the fundamental solution of $Lu=0$ with the
boundary conditions specified before. More developments on this kind of SPDEs and generalized ones can
be found in \cite{Da}, \cite{KZ}, \cite{DM}, \cite{Sa} and the references
therein. \\

Our purpose in this paper is to study some stability, prevalence
and genericity (in the sense of Baire categories) results for the solution of SPDEs
of the form (\ref{integrale}). Specifically, we examine these points :

\begin{enumerate}
\item  Stability of the solution under Lipschitz conditions on the
coefficients

\item  Stability of the solution with respect to the driving process under
pathwise uniqueness

\item  Pathwise uniqueness and convergence of the Picard successive
approximation

\item  Generic properties of the existence, uniqueness and continuous
dependence on the coefficients
\end{enumerate}

\noindent

The paper is organized as follows. In Section 2 we state a H\"{o}lder
regularity of the solution of (\ref{integrale}). Section 3 is devoted to prove some results of
stability of solutions of (\ref{eq0}), first on the coefficients and then
with respect to the driving process. In section 4 we give a necessary and
sufficient conditions which ensure the convergence of the Picard successive
approximation associated to the equation (\ref{integrale}). Section 5 shows
that, the existence, uniqueness and the continuous dependence on the
coefficients are a generic properties in the sense of Baire categories. In
section 6, we study some examples of SPDEs of kind (\ref{eq0}). Finally, an
appendix gathers some technical lemmas which are used through the work. We
always assume that all constants will be denoted by $c$ or $C$ independently
of its value. \noindent In the sequel, we shall refer the equation (\ref{eq1}%
) as Eq($\sigma ,b$). To simplify the notation, we shall write 
\[
b(u)(s,y)=b(s,y,u(s,y))\mbox{ \ and \ }\sigma (u) (s,y)=\sigma
(s,y,u(s,y)).
\]

\section{Definitions and H\"older regularity of the solution}

A weak solution to (\ref{integrale}) is a solution on some filtered space
with respect to some noise $M$, i.e. the noise and space are not specified
in advance. A strong solution to (\ref{integrale}) is a solution which is
adapted with respect to the canonical filtration of the noise $M$.\\

\noindent As solution spaces we consider the spaces defined by:
\begin{definition}
\label{def1} Let $p\geq 2$, a stochastic process $u$ defined on $\Omega
\times \mathbb{R}_{+}\times \mathbb{R}^{d}$, which is jointly measurable and $%
\mathcal{G}_{t}$--adapted, is said to be a solution to the SPDE $(\ref{eq0})$%
, if it is an $\mathbb{R}$--valued fields which satisfies $(\ref{eq1})$ and 
$\sup_{t\in [0,T]}\sup_{x\in \mathbb{R}^{d}}E|u(t,x)|^{p}<+\infty.$
\end{definition}

\noindent Let $f$ be a real valued function $f$ defined on $\mathbb{R}%
_+\times \mathbb{R}^{d}\times \mathbb{R}$.\\
We say that $f$ satisfies (\textbf{L}) if there exists is a constant $c$ such that
\[
\left| f(t,x,r)-f(t,x,v)\right| \leq c|r-v|\mbox{ for
all }t,x,r,v , 
\]
and $f$ satisfies (\textbf{LG}) if 
\[
\left| f(t,x,r)\right| \leq c\left( 1+|r|\right) \mbox{ for all }t,x,r.
\]

\noindent We assume the following set of assumptions: \\

\noindent \textbf{Assumption} \textbf{(A$_{\eta}$)}
\[
\int_{\mathbb{R}^{d}}\frac{\mu(d\xi)}{(1+\vert \xi \vert^2)^\eta}<\infty,
\,\,\, \mbox{for} \,\,\,\eta\in (0,1].
\]
\textbf{Assumption} \textbf{(R.1)}

\begin{enumerate}
\item[(i)]  For any $T>0$, 
\[
\int_{0}^{T}ds\int_{\mathbb{R}^{d}}\mu (d\xi )\left| \mathcal{F}S(s,\cdot )(\xi
)\right| ^{2}<\infty. 
\]

\item[(ii)]  There exist constants $c>0$ and $\delta _{1}>0$ such that for $%
0\leq t_{1}\leq t_{2}\leq T$, 
\[
\sup_{x\in \mathbb{R}^{d}}\int_{t_{1}}^{t_{2}}ds\int_{\mathbb{R}^{d}}\mu (d\xi
)\left| \mathcal{F}S(t_{2}-s,x,\cdot )(\xi )\right| ^{2}\leq c\left|
t_{2}-t_{1}\right| ^{2\delta _{1}}.
\]
\end{enumerate}
\noindent \textbf{Assumption} \textbf{(R.2)}
\begin{enumerate}
\item [ ] For any compact subset $K\subset \mathbb{R}^{d}$ there exist
constants $c>0$ and $\delta _{2}>0$ such that 
\[
\int_{0}^{T}ds\int_{\mathbb{R}^{d}}\mu (d\xi )\left| \mathcal{F}S(s,x+z,\cdot
)(\xi )-\mathcal{F}S(s,x,\cdot )(\xi )\right| ^{2}\leq c\left\| z\right\|
^{2\delta _{2}}\mbox{,} 
\]
for any $x\in \mathbb{R}^{d}$ and $z\in K$.
\end{enumerate}

\noindent It is proved by Dalang \cite{Da} that the assumption $(\mathbf{{A_1%
})}$ together with (\textbf{L}) and (\textbf{LG}) on $\sigma$ and $b$ ensure
the existence and uniqueness of the solution of (\ref{integrale}).

\noindent Let us recall some recent results on the regularity of $u(t,x)$,
which has been proved by Sanz--Sol\'e and Sarr\`a \cite{SaSar} (see also 
\cite{Sa} and \cite{DaSa}). \newline
\noindent Let $\gamma =(\gamma _{1},\gamma _{2})$ such that $\gamma
_{1},\gamma _{2}>0$ and let $K$ be a compact subset of $\mathbb{R}^{d}$. We
denote by $\mathcal{C}^{\gamma }\left( [0,T]\times K;\mathbb{R}\right) $ the
set of $\gamma $--H\"{o}lder continuous functions equipped with the norm
defined by: 
\begin{equation}
\left\| f\right\| _{\gamma ,T,K}=\sup_{(t,x)\in [0,T]\times K}\left|
f(t,x)\right| +\sup_{s\neq t\in [0,T]}\sup_{x\neq y\in K}\frac{\left|
f(t,x)-f(s,y)\right| }{\left| t-s\right| ^{\gamma _{1}}+\left\| x-y\right\|
^{\gamma _{2}}}.
\end{equation}

\begin{theorem}
\label{th0}Assume that $(\mathbf{R.1})$ and $\mathbf{(R.2)}$ hold and let $u$
be a solution to equation $Eq(\sigma ,b)$.

\begin{enumerate}
\item[(i)]  If $b$ and $\sigma $ satisfy $\mathbf{(LG)}$ then $u$ belongs to 
$\mathcal{C}^{\gamma }\left( [0,T]\times K;\mathbb{R}\right) $ a.s. for any $%
\gamma _{i}<\delta _{i}$, $i=1,2$ and for any compact subset $K$ of $\mathbb{R}%
^{d}$.

\item[(ii)]  Moreover $E\left\| u\right\| _{\gamma ,T,K}^{p}<\infty $ for
any $p\geq 2$.
\end{enumerate}
\end{theorem}

\noindent \textbf{Proof.} The proof of (i) follows using Kolmogorov
criterium for more details (see \cite{SaSar}).

\noindent For (ii) we apply the Garcia--Rodemich--Rumsey lemma (see \cite{N}
(1995), p. 237), to obtain 
\[
E(\left\| u\right\| _{\gamma ,T,K}^{p})\leq c_{p,\gamma ,T,d}.
\]
\hfill $\Box$ 

\section{Stability of the solution}

\setcounter{equation}{0} In this section we give some stability results of
the solution of SPDEs of the form (\ref{eq0}). \newline
For a function $f$ defined on $\mathbb{R}_{+}\times \mathbb{R}^{d}\times \mathbb{R}$
and for $T>0$, we set 
\begin{equation}
\left\| f\right\| _{T,\infty }=\sup_{t\in [0,T]}\sup_{x\in \mathbb{R}%
^{d}}\sup_{r\in \mathbb{R}}\left| f(t,x,r)\right| .  \label{norm}
\end{equation}

\subsection{Stability under Lipschitz conditions on the coefficients}

Assume that $(\mathbf{R.1})$ and $\mathbf{(R.2)}$ hold. Let $(\sigma
_{n})_{n\geq 0}$ and $(b_{n})_{n\geq 0}$ be two sequences of functions on $%
\mathbb{R}_{+}\times \mathbb{R}^{d}\times \mathbb{R}$ which satisfy (\textbf{L}) and (%
\textbf{LG}) uniformly in $n$. \newline

\noindent Denote by $\{u_{n}(t,x), t\geq 0, x\in \mathbb{R}^{d}\}$ the unique solution of
equation Eq$(\sigma _{n},b_{n})$ i.e. 
\begin{eqnarray}
u_{n}(t,x) &=&\int_{0}^{t}\int_{\mathbb{R}^{d}}S(t-s,x,y)\sigma
_{n}(u_{n})(s,y)M(ds,dy)  \nonumber  \label{eq2} \\
&&+\int_{0}^{t}ds\int_{\mathbb{R}^{d}}dyS(t-s,x,y)b_{n}(u_{n})(s,y).
\end{eqnarray}
Then, we have the following theorem:

\begin{theorem}
\label{th1}Assume that $(\sigma _{n})_{n\geq 0}$ and $(b_{n})_{n\geq 0}$
converge uniformly respectively to $\sigma $ and $b$ on compact sets of $%
\mathbb{R}_{+}\times \mathbb{R}^{d}\times \mathbb{R}$. Then for any $p\geq 2$ 
\[
\lim_{n\rightarrow +\infty }E\left( \left\| u_{n}-u\right\| _{\gamma
,T,K}^{p}\right) =0, 
\]
where $u$ is the unique solution of $Eq(\sigma ,b).$
\end{theorem}

\noindent The proof of this Theorem is a consequence of the following lemma
and the Theorem \ref{th0}.

\begin{lemma}
\label{lem2} Assume that there exist real valued functions $\sigma $ and $b$
defined on $\mathbb{R}_{+}\times \mathbb{R}^{d}\times \mathbb{R}$ such that 
\begin{equation}
\lim_{n\rightarrow +\infty }\left( \left\| \sigma _{n}-\sigma \right\|
_{T,\infty }+\left\| b_{n}-b\right\| _{T,\infty }\right) =0.  \label{stab}
\end{equation}
Then, for any $p\geq 2$ 
\[
\lim_{n\rightarrow +\infty }\sup_{t\in [0,T]}\sup_{x\in \mathbb{R}%
^{d}}E\left( \left| u_{n}(t,x)-u(t,x)\right| ^{p}\right) =0,
\]
where $u$ is the unique solution of $\mathit{Eq}(\sigma ,b).$
\end{lemma}

\bigskip

\noindent \textbf{Proof.} Without loss of generality assume that $b\equiv 0$%
. For, $p\geq 2$, set 
\[
\varphi_{n}(t,x)=E\left| u_{n}(t,x)-u(t,x)\right|^{p}
\]
and $\phi _{n}(t)=\sup_{x\in \mathbb{R}}\varphi_{n}(t,x)$. Clearly, 
\begin{eqnarray*}
\varphi _{n}(t,x) &\leq &c_{p}E\left| \int_{0}^{t}\int_{\mathbb{R}%
^{d}}S(t-s,x,y)\left( \sigma _{n}(u_{n})(s,y)-\sigma (u_{n})(s,y)\right)
M(ds,dy)\right| ^{p} \\
&&+c_{p}E\left| \int_{0}^{t}\int_{\mathbb{R}^{d}}S(t-s,x,y)\left[ \sigma
(u_{n})(s,y)-\sigma (u)(s,y)\right] M(ds,dy)\right| ^{p}.
\end{eqnarray*}
Then, Burkholder's, H\"{o}lder's inequalities and the property (\textbf{L})
on $\sigma _{n}$ imply 
\[
\phi _{n}(t)\leq c_{p}\left\| \sigma _{n}-\sigma \right\| _{T,\infty
}^{p}\nu _{t}^{\frac{p}{2}}+c_{p}\nu _{t}^{\frac{p}{2}-1}\int_{0}^{t}J(t-s)%
\phi _{n}(s)ds, 
\]
where 
\begin{equation}  \label{Jmu}
J(s)=\int_{\mathbb{R}^{d}}\mu (d\xi )\left| \mathcal{F}S(s,\cdot )(\xi )\right|
^{2}\mbox{ \ and \ }\nu _{t}=\int_{0}^{t}J(t-s)ds.
\end{equation}
Hence 
\[
\phi _{n}(t)\leq c_{p,T}\left( \left\| \sigma _{n}-\sigma \right\|
_{T,\infty }^{p}+\int_{0}^{t}J(t-s)\phi _{n}(s)ds\right) . 
\]
Therefore the hypothesis \textbf{(R.1)}-(i) and the lemma \ref{lemGG} yield 
\[
\sup_{t\in [0,T]}\phi _{n}(t)\leq c_{p,T}\left\| \sigma _{n}-\sigma \right\|
_{T,\infty }^{p}. 
\]
Now, apply (\ref{stab}) to complete the proof of the lemma. \hfill $\Box$

\bigskip

\noindent \textbf{Proof of Theorem \ref{th1}.} It suffices to prove that
the sequence $u_{n}-u$ satisfies the properties (\textbf{P}$_{1}$) and (%
\textbf{P}$_{2}$) of lemma \ref{lem1}. Clearly, by Theorem \ref{th0}, $%
u_{n}-u$ satisfy the property (\textbf{P}$_{1}$) of Lemma \ref{lem1}. The
property (\textbf{P}$_{2}$) is given by Lemma \ref{lem2}. Therefore the
proof of Theorem \ref{th1} follows from the above properties. \hfill $\Box$

\begin{definition}
\label{def2}We say that the pathwise uniqueness property $\mathbf{(PU)}$ holds for
equation (\ref{eq1}) if whenever $(u,M,(\Omega ,\mathcal{G},P),\mathcal{G}%
_{t})$ and $(u^{\prime },M^{\prime },(\Omega ,\mathcal{G},P),\mathcal{G}%
_{t}^{\prime })$ are two weak solutions of equation (\ref{eq1}) such that $M
\equiv M^{\prime }$\ $P$--a.s., then $u\equiv u^{\prime }$ $P$--a.s.
\end{definition}

\noindent In the next, we state a variant of the Theorem \ref{th1}. Let us
consider a family of functions depending on a parameter $\lambda \in \mathbb{R}$%
, and consider the stochastic partial differential equation: 
\begin{eqnarray}
u^{\lambda }(t,x) &=&\varphi \left( \lambda \right)+\int_{0}^{t}\int_{\mathbb{R}%
^{d}}S(t-s,x,y)\sigma _{\lambda }( u^{\lambda }) (s,y)M(ds,dy)  \nonumber
\label{eqpara} \\
&&+\int_{0}^{t}ds\int_{\mathbb{R}^{d}}dyS(t-s,x,y)b_{\lambda }( u^{\lambda })
(s,y), \\
u^{\lambda }(0,x) &=&\varphi \left( \lambda \right) \mbox{ for }x\in
\partial (\mathbb{R}^{d}), \;\, \mbox{where $\varphi$ is a given function.} 
\nonumber
\end{eqnarray}
Then, we have the following theorem:

\begin{theorem}
\label{th1p} Suppose that $\sigma _{\lambda }(t,x,r) $ and $b_{\lambda
}(t,x,r) $ are continuous with respect to their arguments. Further, suppose
that $\varphi$ is continuous at $\lambda _{0}\in \mathbb{R},$ and for each $T>0$
and each compact subset $K$ of $\mathbb{R}^{d}$ there exists a constant $c>0$
such that for all $r\in \mathbb{R}$ 
\begin{eqnarray*}
\sup_{\lambda \in \mathbb{R}}\sup_{t\in [0,T]}\sup_{x\in \mathbb{R}^{d}}\left(
\left| \sigma _{\lambda}(t,x,r) \right| +\left| b_{\lambda}(t,x,r) \right|
\right) &\leq & c\left( 1+\left| r\right| \right) 
\end{eqnarray*}
and
\begin{eqnarray*}
\lim_{\lambda \rightarrow \lambda _{0}}\left( \left\| \sigma _{\lambda
}-\sigma _{\lambda _{0}}\right\| _{T,\infty }+\left\| b_{\lambda
}-b_{\lambda _{0}}\right\| _{T,\infty }\right) &=&0,
\end{eqnarray*}
where $\Vert \cdot \Vert _{T,\infty }$ has been defined in $(\ref{norm})$.
Then, under $\mathbf{(PU)}$ for the equation $(\ref{eqpara})$ at $\lambda
_{0}$ we have: 
\[
\lim_{\lambda \rightarrow \lambda _{0}}E\left[ \left\| u^{\lambda
}-u^{\lambda _{0}}\right\| _{T,\infty }^{2}\right] =0\,\ \mbox{for every}%
\;T\geq 0. 
\]
\end{theorem}

\bigskip

\noindent \textbf{Proof.} Similar to the proof of Theorem \ref{th1}.\hfill $%
\Box$

\subsection{Stability with respect to the driving process under \textbf{(PU)}}

In this subsection, we consider SPDEs driven by spatially correlated noise.
We prove a continuity result with respect to the driving processes, when the
pathwise uniqueness of solutions holds.

\noindent Let $\{M^{n}\}_{n\geq 0}$ be a sequence of continuous $(\mathcal{G}%
_{t},P)$--martingale measure, with $M^0=M$ and $\sigma ,b:\mathbb{R}_{+}\times \mathbb{R}^{d}\times \mathbb{R}
\longrightarrow \mathbb{R}$ be continuous functions satisfying (\textbf{LG}). Define the sequence
\begin{eqnarray*}
u^{n}(t,x) &=&\int_{0}^{t}\int_{\mathbb{R}^{d}}S(t-s,x,y)\sigma \left(
u^{n}\right) (s,y)M^{n}(ds,dy) \\
&&+\int_{0}^{t}ds\int_{\mathbb{R}^{d}}dyS(t-s,x,y)b\left( u^{n}\right) (s,y).
\end{eqnarray*}

\noindent Suppose that $\{M^{n}\}_{n\geq 0}$ satisfy the following
conditions:

\begin{quote}
$(\mathbf{H.1})$ The family $\{M^{n}\}_{n\geq 0}$ is bounded in probability
in $\mathcal{C}\left( [0,T]\times K\right) .$

$(\mathbf{H.2})$ $M^{n}-M^{0}\longrightarrow _{n\rightarrow +\infty }0$ in
probability on $\mathcal{C}\left( [0,T]\times K\right).$
\end{quote}

\noindent Then, we have the following theorem:

\begin{theorem}
\label{thsn2} Suppose that $\left( \mathbf{R.1}\right) $, $\left( \mathbf{R.2}\right) $%
, $\left( \mathbf{H.1}\right) $, $\left( \mathbf{H.2}\right) $ and $\mathbf{%
(PU)}$ hold and that $Eq(\sigma ,b)$ is non--degenerate. Assume,
moreover, that $\frac{\partial \sigma }{\partial r}$ is a locally bounded
functions of $(t,x,r)$ and that it is Lipschitz continuous in $r\in \mathbb{R}$%
. Then for any $\varepsilon >0$, 
\[
\lim_{n\rightarrow +\infty }P\left( \left\| u^{n}-u\right\| _{\gamma
,T,K}>\varepsilon \right) =0.
\]
\end{theorem}

\noindent The main tool used in the proofs is the Skorokhod representation
theorem given by the following:

\begin{lemma}
\label{lemtc1}$([11]\mbox{ p. 9})$ Let $\left( \mathcal{X},\rho \right) $ be
a complete separable metric space, $\{P_{n}:n\geq 1\}$ and $P$ be
probability measures on $\left( \mathcal{X},\mathcal{B}\left( \mathcal{X}%
\right) \right) $ such that $P_{n}\rightarrow _{n\rightarrow +\infty }P$.
Then, on a probability space $(\widehat{\Omega },\widehat{\mathcal{G}},%
\widehat{P})$, we can construct $\mathcal{X}$--valued random variables $%
\{u_{n}:n\geq 1\}$ and $u$ such that:

$(i)$ $P_{n}=\widehat{P}_{u_{n}}$, $n=1,2,...$ and $P=\widehat{P}_{u}.$

$(ii)$ $u_{n}$ converges to $u$ $\widehat{P}$--a.s.
\end{lemma}

We will make use of the following result, which gives us a criteria for the
tightness of sequences of laws associated to continuous processes.

\begin{lemma}
\label{lemtc2}$([11]\mbox{ p. 18})$ Let $\{u_{n}(t,x):n\geq 1\}$ be a
sequence of real valued continuous processes satisfying the following two
conditions:

$(i)$ There exist positive constants $C$ and $q$ such that 
\[
\sup_{n\geq 1}E\left[ \left| u_{n}(0,x_0)\right| ^{q}\right] \leq C \; 
\mbox{for
some given } \; x_0.
\]

$(ii)$ There exist positive constants $p$, $\beta _{1}$, $\beta _{2}$, $%
C_{T} $ such that: 
\[
\sup_{n\geq 1}E\left[ \left| u_{n}(t,x)-u_{n}(s,y)\right| ^{p}\right] \leq
C_{T}\left( \left| t-s\right| ^{1+\beta _{1}}+\left\| x-y\right\| ^{d+\beta
_{2}}\right) 
\]
for every $s,t\in [0,T]$ and $x,y\in \mathbb{R}^{d}$.

Then, there exist a subsequence $(n_{k})_{k\geq 1}$, a probability space $(%
\widehat{\Omega },\widehat{\mathcal{G}},\widehat{P})$ and real valued
continuous processes $\widehat{u}_{n_{k}}$, $k=1,2,\ldots $ and $\widehat{u}$
defined on it such that:

\begin{enumerate}
\item  The two random field $\widehat{u}_{n_{k}}$ and $u_{n_{k}}$ have the
same law.

\item  $\widehat{u}_{n_{k}}(t,x)$ converge to $\widehat{u}(t,x)$ uniformly
on every compact subset of $\mathbb{R}_{+}\times \mathbb{R}^{d}$ $\widehat{P}$%
--a.s.
\end{enumerate}
\end{lemma}
The following lemma is a variant of the Lemma 4.3 in (\cite{GN} p. 744--745) in which we replace the Brownian sheet with a martingale measure.
\begin{lemma}
\label{lemst}For every $n\geq 0$ let $\{z^{n}(t,x):t\in 
\mathbb{R}_{+}\times \mathbb{R}^{d}\}$ be a family of continuous $\mathcal{G}%
_{t}^{n}$--adapted random field and let $M^{n}$ be a martingale measure
carried by some filtered probability space $(\Omega ,\mathcal{G},\mathcal{G}%
_{t}^{n},P)$. Assume that for every $\varepsilon >0$, $T>0$ and $K$ a
compact subset of $\mathbb{R}^{d}$: 
\[
\lim_{n\rightarrow +\infty }P\left( \left\| z^{n}-z^{0}\right\| _{\gamma
,T,K}+\left\| M^{n}-M^{0}\right\| _{T,\infty }>\varepsilon \right) =0. 
\]
Let $h(t,x,r)$ be a bounded Borel function of $(t,x,r)\in \mathbb{R}_{+}\times 
\mathbb{R}^{d}\times \mathbb{R}$. \newline
Then the following assertions hold:

$(i)$ If $h$ is continuous in $r\in \mathbb{R}$, then 
\begin{equation}
\lim_{n\rightarrow +\infty }I_{n}(t,x)=I_{0}(t,x)\mbox{ \ and \ }%
\lim_{n\rightarrow +\infty }J_{n}(t,x)=J_{0}(t,x)  \label{kry}
\end{equation}
in probability for every $t\geq 0$ and every $x\in \mathbb{R}^{d}$, where 
\[
I_{n}(t,x)=\int_{0}^{t}\int_{\mathbb{R}^{d}}S(t-s,x,y)h(s,y,z^{n}(s,y)dyds 
\]
and 
\[
J_{n}(t,x)=\int_{0}^{t}\int_{\mathbb{R}%
^{d}}S(t-s,x,y)h(s,y,z^{n}(s,y)M^{n}(ds,dy) 
\]

$(ii)$ If for almost every $(t,x)\in \mathbb{R}_{+}\times \mathbb{R}^{d}$ the law $%
Q_{t,x}^{n}$ of $z^{n}(t,x)$ is absolutely continuous with respect to the
Lebesgue measure $\lambda $ on $\mathbb{R}$ and the density $p_{t,x}^{n}=\frac{%
dQ_{t,x}^{n}}{d\lambda }$ satisfies for some $\alpha >1$%
\[
\sup_{n\geq 0}\int_{0}^{T}\int_{\mathbb{R}^{d}}\int_{\mathbb{R}}\left(
p_{t,x}^{n}\left( r\right) \right) ^{\alpha }drdxdt<\infty . 
\]
Then $(\ref{kry})$ also hold in probability for every $t\geq 0$ and every $%
x\in \mathbb{R}^{d}.$
\end{lemma}

\noindent \textbf{Proof of Theorem \ref{thsn2}.} Suppose that the
conclusion of our theorem is false. Then there exists $\varepsilon >0$ and a
subsequence $\left( n_{k}\right) _{k\geq 0}$ such that 
\[
\inf_{n_{k}}P\left( \left\| u^{n_{k}}-u\right\| _{\gamma ,T,K}>\varepsilon
\right) \geq \varepsilon . 
\]
Clearly the family $Z^{n}=\left( u^{n},u,M^{n},M\right) $ is tight in $(%
\mathcal{C}([0,T]\times \mathbb{R}^{d},\mathbb{R}))^{4}$. Then, by Skorokhod's
representation theorem, there exist a probability space $(\widehat{\Omega },%
\widehat{\mathcal{G}},\widehat{P})$ and $\widehat{Z}^{n_{k}}=(\widehat{u}%
^{n},\widetilde{u}^{n},\widehat{M}^{n},\widetilde{M}^{n})$ which satisfy:

\begin{quote}
$i)$ $Law(Z^{n_{k}})=Law(\widehat{Z}^{n_{k}})$

$ii)$ There exists a subsequence $(\widehat{Z}^{n_{k}})_{k}$ also denoted by 
$(\widehat{Z}^{n})_{n}$ which converges $\widehat{P}$--a.s. in $[\mathcal{C}%
([0,T]\times \mathbb{R}^{d},\mathbb{R})]^{4}$ to $\widehat{Z}=(\widehat{u},%
\widetilde{u},\widehat{M},\widetilde{M}).$
\end{quote}

\noindent Let $\mathcal{G}_{t}^{n}$ denotes the completion of the $\sigma $%
--algebra generated by $\{\widehat{Z}_{s}^{n}: s\le t\}$ and set $\widehat{%
\mathcal{G}}_{t}^{n}=\cap _{s>t}\mathcal{G}_{s}^{n}.$\newline

\noindent In an analogous manner we define the $\sigma $--algebra $\{%
\widehat{\mathcal{F}}_{t}:t\in [0,T]\}$ for the limiting process $%
\widehat{Z}$. Then $(\widehat{\Omega },\widehat{\mathcal{G}},\widehat{%
\mathcal{G}}_{t}^{n},\widehat{P})$ (resp.$\ (\widehat{\Omega },\widehat{%
\mathcal{G}},\widehat{\mathcal{G}}_{t},\widehat{P})$) are stochastic basis
and $\widehat{M}^{n}$, $\widetilde{M}^{n}$ (resp. $\widehat{M}$, $\widetilde{%
M}$) are $\widehat{\mathcal{G}}_{t}^{n}$ (resp. $\widehat{\mathcal{G}}_{t}$%
)--continuous martingale measures. Moreover the two random fields $\widehat{u%
}^{n}$ and $\widetilde{u}^{n}$ satisfy the following SPDEs: 
\begin{eqnarray*}
\widehat{u}^{n}(t,x) &=&\int_{0}^{t}\int_{\mathbb{R}^{d}}S(t-s,x,y)\sigma
\left( \widehat{u}^{n}\right) (s,y)\widehat{M}^{n}(ds,dy) \\
&&+\int_{0}^{t}ds\int_{\mathbb{R}^{d}}dyS(t-s,x,y)b\left( \widehat{u}%
^{n}\right) (s,y).
\end{eqnarray*}
and 
\begin{eqnarray*}
\widetilde{u}^{n}(t,x) &=&\int_{0}^{t}\int_{\mathbb{R}^{d}}S(t-s,x,y)\sigma
\left( \widetilde{u}^{n}\right) (s,y)\widetilde{M}^{n}(ds,dy) \\
&&+\int_{0}^{t}ds\int_{\mathbb{R}^{d}}dyS(t-s,x,y)b\left( \widetilde{u}%
^{n}\right) (s,y)
\end{eqnarray*}
on the same stochastic basis $(\widehat{\Omega },\widehat{\mathcal{G}},%
\widehat{\mathcal{G}}_{t}^{n},\widehat{P}).$

\noindent By using the Lemma \ref{lemst}, we see that the limiting processes
satisfy the following equations: 
\begin{eqnarray*}
\widehat{u}(t,x) &=&\int_{0}^{t}\int_{\mathbb{R}^{d}}S(t-s,x,y)\sigma \left( 
\widehat{u}\right) (s,y)\widehat{M}(ds,dy) \\
&&+\int_{0}^{t}ds\int_{\mathbb{R}^{d}}dyS(t-s,x,y)b\left( \widehat{u}\right)
(s,y).
\end{eqnarray*}
and 
\begin{eqnarray*}
\widetilde{u}(t,x) &=&\int_{0}^{t}\int_{\mathbb{R}^{d}}S(t-s,x,y)\sigma \left( 
\widetilde{u}\right) (s,y)\widetilde{M}(ds,dy) \\
&&+\int_{0}^{t}ds\int_{\mathbb{R}^{d}}dyS(t-s,x,y)b\left( \widetilde{u}\right)
(s,y).
\end{eqnarray*}
Due to $\left( \mathbf{H.2}\right) $ it is easy to see that $\widehat{M}=%
\widetilde{M}$, $\widehat{P}$--a.s.\newline

\noindent Hence by the pathwise uniqueness, $\widehat{u}$ and $\widetilde{u}$
are indistinguishable. This contradicts our assumption. Therefore $u^{n}$
converges to the unique solution $u$. \hfill $\Box$

\section{Pathwise uniqueness and Picard's successive approximation}

Let $\sigma $ and $b$ satisfy \textbf{(LG)} and are continuous. We consider
the SPDEs (\ref{eq1}). The sequence of the Picard successive approximation associated
to (\ref{eq1}) is defined as follows: 
\begin{eqnarray} \label{eqip}
\left\{ 
\begin{array}{ccc}
u^{0}=0 && \\ 
u^{n+1}(t,x)&=&\dint_{0}^{t}\dint_{\mathbb{R}^{d}}S(t-s,x,y)\sigma \left(
u^{n}\right) (s,y)M(ds,dy) \\ 
&&\\ && +\dint_{0}^{t}\dint_{\mathbb{R}^{d}}S(t-s,x,y)b\left(
u^{n}\right) (s,y)dyds.
\end{array}
\right. 
\end{eqnarray}
If we assume that the coefficients $\sigma $ and $b$ satisfy the condition 
\textbf{(L)}, then the sequence $(u^{n})_{n\geq 0}$ converges in $%
L^{p}(\Omega )$ (as $n\rightarrow \infty $) and gives an effective way for
the construction of the unique solution $u$ of equation (\ref{eq1}) (see for
instance \cite{Da}). \\

\noindent Now, if we drop the Lipschitz condition on the coefficients and assume only
that equation (\ref{eq1}) admits a unique strong solution, does the sequence 
$(u^{n})_{n\geq 0}$ converge to $u$ ? The answer is negative even in the
deterministic case, (see \cite{D} p. 114--124).\\

\noindent The aim of the following theorem is to establish an additional necessary and
sufficient condition which ensures the convergence of the Picard successive
approximation.

\begin{theorem}
\label{thip}Let $\sigma $ and $b$ be continuous functions satisfying $%
\mathbf{(LG)}$. Assume further that $\mathbf{(R.1)}$, $\mathbf{(R.2)}$ and $%
\mathbf{(PU)}$ hold for the equation $(\ref{eq1})$. Then $\left(
u^{n}\right) _{n\geq 0}$ converges in $L^{p}(\Omega ;\mathcal{C}^{\gamma
}([0,T]\times K,\mathbb{R}))$, $(p\geq 2, \gamma =(\gamma _{1},\gamma _{2})$
with $\gamma _{1}<\delta _{1}$ and $\gamma _{2}<\delta _{2})$ to the unique
solution of $(\ref{eq1})$ if and only if $\left\| u^{n+1}-u^{n}\right\|
_{\gamma ,T,K}$ converges to $0$, as $n \rightarrow \infty$, in $%
L^{p}(\Omega )$.
\end{theorem}

First, we show tightness of the sequence $u^{n}$.

\begin{lemma}
\label{lemip}Let $\left( u^{n}\right) _{n\geq 0}$ be defined by $(\ref{eqip})$. Then $u^{n}$ is tight in $\mathcal{C}([0,T]\times \mathbb{R}^{d},\mathbb{R})$.
Moreover $\sup_{n\geq 0}E[\left\| u^{n}\right\| _{\gamma ,T,K}^{p}]<+\infty,$
for every $p\geq 2$ and any $\gamma =(\gamma _{1},\gamma _{2})$ such that $%
\gamma _{1}<\delta _{1}$ and $\gamma _{2}<\delta _{2}$.
\end{lemma}

\noindent \textbf{Proof.} For all $t>0$ and $n>1$, we have 
\begin{eqnarray*}
\left| u^{n}(t,x)\right| ^{p} &\leq &c_{p}\left| \int_{0}^{t}\int_{\mathbb{R}%
^{d}}S(t-s,x,y)\sigma( u^{n-1}) (s,y)M(ds,dy)\right| ^{p} \\
&&+c_{p}\left| \int_{0}^{t}ds\int_{\mathbb{R}^{d}}dyS(t-s,x,y)b( u^{n-1})
(s,y)\right| ^{p}.
\end{eqnarray*}
Burkholder and H\"{o}lder inequalities provide the following estimate 
\[
E\left| u^{n}(t,x)\right| ^{p}\leq c_{p}\nu _{t}^{\frac{p}{2}%
-1}\int_{0}^{t}\left( 1+\sup_{y\in \mathbb{R}^{d}}E\left| u^{n-1}(s,y)\right|
^{p}\right) J(t-s)ds. 
\]
Set $\phi _{n}(t)=\sup_{x\in \mathbb{R}^{d}}E\left| u^{n}(t,x)\right| ^{p}$.
Therefore 
\begin{eqnarray*}
\phi _{n}(t) &\leq &c_{p,T}\int_{0}^{t}\left( 1+\phi _{n-1}(s)\right)
J(t-s)ds \\
&\leq &c_{p,T}+c_{p,T}\int_{0}^{t}\phi _{n-1}(s)J(t-s)ds.
\end{eqnarray*}
By \textbf{(R.1)}-(i) and the Lemma 15. p. 22 in Dalang (1999), we deduce
that 
\begin{equation}
\sup_{n\geq 0}\sup_{t\in [0,T]}\sup_{x\in \mathbb{R}^{d}}E\left|
u^{n}(t,x)\right| ^{p}=\sup_{n\geq 0}\sup_{t\in [0,T]}\sup_{x\in \mathbb{R}%
^{d}}\phi _{n}(t)\leq c_{p,T}.  \label{bi}
\end{equation}
To prove the tightness, we write 
\begin{eqnarray*}
u^{n}(t_{2},x_{2})-u^{n}(t_{1},x_{1}) &=&\int_{0}^{t_{2}}\int_{\mathbb{R}%
^{d}}S(t_{2}-s,x_{2},y)\sigma (u^{n-1})(s,y)M(ds,dy) \\
&&-\int_{0}^{t_{1}}\int_{\mathbb{R}^{d}}S(t_{1}-s,x_{1},y)\sigma
(u^{n-1})(s,y)M(ds,dy) \\
&=&\int_{0}^{t_{1}}\int_{\mathbb{R}^{d}}\Lambda
(t_{1},t_{2},x_{1},x_{2},s,y)\sigma (u^{n-1})(s,y)M(ds,dy) \\
&&+\int_{t_{1}}^{t_{2}}\int_{\mathbb{R}^{d}}S(t_{2}-s,x_{2},y)\sigma
(u^{n-1})(s,y)M(ds,dy).
\end{eqnarray*}
Taking the expectation, we have 
\begin{eqnarray*}
&& \lefteqn{E\left| u^{n}(t_{2},x_{2})-u^{n}(t_{1},x_{1})\right| ^{p}} \\
& \leq & c_{p}E\left| \int_{0}^{t_{1}}\int_{\mathbb{R}^{d}}\Lambda
(t_{1},t_{2},x_{1},x_{2},s,y)\sigma (u^{n-1})(s,y)M(ds,dy)\right| ^{p} \\
&& + c_{p}E\left| \int_{t_{1}}^{t_{2}}\int_{\mathbb{R}%
^{d}}S(t_{2}-s,x_{2},y)\sigma (u^{n-1})(s,y)M(ds,dy)\right| ^{p} \\
& \leq & c_{p}\left[\left( \int_{0}^{t_{1}}ds\int_{\mathbb{R}^{d}}\mu
(d\xi )\left| \mathcal{F}\Lambda (t_{1},t_{2},x_{1},x_{2},s,\cdot )(\xi
)\right| ^{2}\right) ^{\frac{p}{2}-1} \right. \\
&&\left.\,\,\,\,\,\,\,\,\times \int_{0}^{t_{1}}ds\left( 1+\phi
_{n-1}(s)\right) \int_{\mathbb{R}^{d}}\mu (d\xi )\left| \mathcal{F}\Lambda
(t_{1},t_{2},x_{1},x_{2},s,\cdot )(\xi )\right| ^{2}\right] \\
&& + c_{p}\left[\left( \int_{t_{1}}^{t_{2}}ds\int_{\mathbb{R}^{d}}\mu
(d\xi )\left| \mathcal{F}S(t_{2}-s,x_{2},\cdot )(\xi )\right| ^{2}\right) ^{%
\frac{p}{2}-1} \right. \\
&&\left.\,\,\,\,\,\,\,\,\times \int_{t_{1}}^{t_{2}}ds\left( 1+\phi
_{n-1}(s)\right) \int_{\mathbb{R}^{d}}\mu (d\xi )\left| \mathcal{F}%
S(t_{2}-s,x_{2},\cdot )(\xi )\right| ^{2}\right].
\end{eqnarray*}
Hence, by (\ref{bi}), $(\mathbf{R.1})$ and $(\mathbf{R.2)}$, we have
\[
E\left| u^{n}(t_{2},x_{2})-u^{n}(t_{1},x_{1})\right| ^{p}\leq c_{p}\left(
\left| t_{2}-t_{1}\right| ^{\delta _{1}p}+\left\| x_{2}-x_{1}\right\|
^{\delta _{2}p}\right) . 
\]
Therefore $u^n$ is tight in $\mathcal{C}([0,T]\times \mathbb{R}^{d},\mathbb{R})$.%
\newline
Now, by Kolmogorov criterium we deduce that $\sup_{n\geq 0}E[\left\|
u^{n}\right\| _{\gamma ,T,K}^{p}]<+\infty $, (see \cite{N}, p. 238).\hfill $%
\Box$ \\

\noindent \textbf{Proof of Theorem \ref{thip}.} Suppose that $\left\|
u^{n+1}-u^{n}\right\| _{\gamma ,T,K}$ converges to $0$, as $n \rightarrow
\infty$, in $L^{p}(\Omega )$, $(p\geq 2)$ and that there is some $%
\varepsilon >0$, and a sequence $(n_{k})_k$ such that: 
\[
\inf_{n_{k}}E\left( \left\| u^{n_{k}}-u\right\| _{\gamma ,T,K}^{p}\right)
\geq \varepsilon . 
\]
According to Lemma \ref{lemip}, the family $\left( u^{n},u^{n+1},u,M\right) $
satisfies conditions $(i)$ and $(ii)$ of Lemma \ref{lemtc2}. Then by the
Skorokhod selection theorem, there exists some probability space $(\widehat{%
\Omega },\widehat{\mathcal{G}},\widehat{P})$ carrying a sequence of
stochastic processes $(\widehat{u}^{n},\widetilde{u}^{n+1},\overline{u}^{n},%
\widehat{M}^{n})$, with the following properties:

\begin{quote}
P.1. For each $n\in \mathbb{N}$, the two random field $(\widehat{u}^{n},%
\widetilde{u}^{n+1},\overline{u}^{n},\widehat{M}^{n})$ and $\left(
u^{n},u^{n+1},u,M\right) $ have the same law for each $n\in \mathbb{N}.$

P.2. There exists a subsequence $(n_{k})_{k\geq 0}$ such that $(\widehat{u}%
^{n_{k}},\widetilde{u}^{n_{k}+1},\overline{u}^{n_{k}},\widehat{M}^{n_{k}})$
converges to $(\widehat{u},\widetilde{u},\overline{u},\widehat{M}) $
uniformly on every compact subset of $\mathbb{R}_{+}\times \mathbb{R}^{d}$ $%
\widehat{P}$--a.s.
\end{quote}

But we know that $u^{n+1}-u^{n}$ converges to $0$, then we can show easily
that $\widehat{u}=\widetilde{u}$, $\widehat{P}$--a.s. If we denote 
\[
\widehat{\mathcal{G}}_{t}^{n}=\sigma \left( \widehat{u}^{n}(s,y),\;\overline{%
u}^{n}(s,y),\;\widehat{M}^{n}(s,y)\ ;\ s\leq t,\;y\in K\right) 
\]
and 
\[
\widehat{\mathcal{G}}_{t}=\sigma \left( \widehat{u}(s,y),\;\overline{u}%
(s,y),\;\widehat{M}(s,y)\ ;\ s\leq t,\;y\in K\right) , 
\]
then $(\widehat{M}^{n},\widehat{\mathcal{G}}_{t}^{n})$ and $(\widehat{M},%
\widehat{\mathcal{G}}_{t})$ are gaussian processes (even martingales
measures) which have the same law as $M$.\newline

According to the property P.1. and the fact that $u^{n}$ and $u$ satisfy
respectively (\ref{eqip}) and (\ref{eq1}) with the same initial condition $%
u_{0}\equiv 0$, it can be proved following the method used by Krylov \cite
{Kr} for diffusions p. 89, that $\forall \;n\in \mathbb{N}$, $\forall \;t\geq 0$
and $x\in \mathbb{R}^{d}$%
\begin{eqnarray*}
&&E\left| \widehat{u}^{n}(t,x)-\int_{0}^{t}\int_{\mathbb{R}^{d}}S(t-s,x,y)%
\sigma \left( \widehat{u}^{n}\right) (s,y)\widehat{M}^{n}(ds,dy)\right. \\
&&\;\;\;\;\;\;\;\;\;\;\;\;\;\;\;\;\;\;\;\left. -\int_{0}^{t}ds\int_{\mathbb{R}%
^{d}}dyS(t-s,x,y)b\left( \widehat{u}^{n}\right) (s,y)\right| ^{2}=0.
\end{eqnarray*}
In other words, $\widehat{u}^{n}$ satisfies the stochastic integral equation
: 
\begin{eqnarray}
\widehat{u}^{n}(t,x) &=&\int_{0}^{t}\int_{\mathbb{R}^{d}}S(t-s,x,y)\sigma
\left( \widehat{u}^{n}\right) (s,y)\widehat{M}^{n}(ds,dy)  \label{eq1n} \\
&&+\int_{0}^{t}ds\int_{\mathbb{R}^{d}}dyS(t-s,x,y)b\left( \widehat{u}%
^{n}\right) (s,y).  \nonumber
\end{eqnarray}
Similarly $\overline{u}^{n}$ satisfies the equation (\ref{eq1n}) that is 
\begin{eqnarray*}
\overline{u}^{n}(t,x) &=&\int_{0}^{t}\int_{\mathbb{R}^{d}}S(t-s,x,y)\sigma
\left( \overline{u}^{n}\right) (s,y)\widehat{M}^{n}(ds,dy) \\
&&+\int_{0}^{t}ds\int_{\mathbb{R}^{d}}dyS(t-s,x,y)b\left( \overline{u}%
^{n}\right) (s,y).
\end{eqnarray*}
By using the property P.2. and a limit theorem of Skorokhod \cite{S} p. 32,
it holds that 
\[
\int_{0}^{t}\int_{\mathbb{R}^{d}}S(t-s,x,y)\sigma \left( \widehat{u}%
^{n_{k}}\right) (s,y)\widehat{M}^{n_{k}}(ds,dy)
\]
and 
\[
\int_{0}^{t}ds\int_{\mathbb{R}^{d}}dyS(t-s,x,y)b\left( \widehat{u}%
^{n_{k}}\right) (s,y)
\]
converge, respectively, in probability (as $k\rightarrow \infty$) to

\[
\int_{0}^{t}\int_{\mathbb{R}^{d}}S(t-s,x,y)\sigma \left( \widehat{u}\right)
(s,y)\widehat{M}(ds,dy) 
\]
and 
\[
\int_{0}^{t}ds\int_{\mathbb{R}^{d}}dyS(t-s,x,y)b\left( \widehat{u}\right)
(s,y). 
\]

Therefore $\widehat{u}$ and $\overline{u}$ satisfy the same stochastic
partial differential equation (\ref{eq1}), on the space $(\widehat{\Omega },%
\widehat{\mathcal{G}},\widehat{P})$, with the same gaussian noise $\widehat{M%
}$ and the same initial condition. Then, by the pathwise uniqueness
property, we conclude that $\widehat{u}(t,x)=\overline{u}(t,x)$, $\forall
\;t,x\,\widehat{P}$--a.s.\newline
By the uniform integrability, it holds that for some $\varepsilon >0$: 
\begin{eqnarray*}
\varepsilon &\leq &\liminf_{n\in \mathbb{N}}E(\left\| u^{n}-u\right\| _{\gamma
,T,K}^{p}) \\
&\leq &\liminf_{k\in \mathbb{N}}\widehat{E}\left( \left\| \widehat{u}^{n_{k}}-%
\overline{u}^{n_{k}}\right\| _{\gamma ,T,K}^{p}\right) =\widehat{E}\left( \left\| \widehat{u}-\overline{u}\right\| _{\gamma
,T,K}^{p}\right) ,
\end{eqnarray*}
which is a contradiction.\hfill $\Box$

\begin{remark}
Note that under $(\mathbf{R.1})$, $(\mathbf{R.2)}$ and $(\mathbf{PU})$, the
series $\sum_{n\geq 0}\left( u^{n+1}-u^{n}\right) $ converges in $L^{p}(\Omega ;\mathcal{C}^{\gamma
}([0,T]\times K,\mathbb{R}))$, $(p\geq 2)$ if and only
if $\left( u^{n+1}-u^{n}\right) _{n\geq 0}$ converges to $0$ in $L^{p}(\Omega ;\mathcal{C}^{\gamma
}([0,T]\times K,\mathbb{R}))$.
\end{remark}

The main result of this section is

\begin{theorem}
\label{thips}Let $\sigma $ and $b$ be continuous functions satisfying $%
\mathbf{(LG)}$. Let $\left( u^{n}\right) _{n\geq 0}$\ be given by $(\ref
{eqip})$. Assume further that $(\mathbf{R.1})$, $(\mathbf{R.2)}$ and $(\mathbf{PU})$ holds for the equation $(\ref
{eq1})$, then for any $p\geq 2$ 
\[
\lim_{n\rightarrow +\infty }E\left[ \left\| u^{n}-u\right\| _{\gamma
,T,K}^{p}\right] =0, 
\]
where $u$ is the unique solution of $(\ref{eq1})$.
\end{theorem}

\noindent \textbf{Proof of Theorem \ref{thips}.} Similar than the proof of
Theorem \ref{thip} so it is omitted.\hfill $\Box$

\section{Genericity of the existence and uniqueness}

As we have seen in previous sections, the pathwise uniqueness property plays
a key role in the proof of many stability results. It is then quite natural
to raise the question whether the set of all ``\textit{nice}'' functions ($%
\sigma ,b)$ for which the pathwise uniqueness holds for the stochastic
partial differential equation Eq($\sigma ,b$) is larger than its complement,
in a sense to be specified. To make the question meaningful let us recall
what we mean by the generic property.

\noindent A property $\mathcal{P}$ is said to be generic for a class of
stochastic partial differential equations $\mathcal{E}$, if $\mathcal{P}$ is
satisfied by each equation in $\mathcal{E}\setminus \mathcal{A}$, where $%
\mathcal{A}$ is a set of first category (in the sense of Baire) in $\mathcal{%
E}$. Results on generic properties for ordinary differential equations seems
to go back to an old paper of Orlicz \cite{O}, see also Lasota and Yoke \cite
{LY}. The investigation of such questions for stochastic differential
equations is carried out in Heunis \cite{H} and Alibert and Bahlali \cite{AB}. In this section, we show that the subset of continuous and bounded
coefficients for which existence, pathwise uniqueness holds for equation Eq($\sigma
,b$) is a residual set. The proof is based essentially on Theorem \ref{th1p}. Moreover it does not use the oscillation function introduced by Lasota and
Yorke in \cite{LY} in ordinary differential equations and used in stochastic
differential equations by Heunis \cite{H} and in farther development for the
generic property of stochastic differential equations by Bahlali \emph{et al}. \cite{BMO1}. See also
Bahlali \emph{et al}. \cite{BMO2} for backward stochastic differential
equations.\newline

\noindent In this section we improve the result obtained in \cite{BEE}  by
considering small spaces which contain the solutions of SPDEs in any space
dimension but driven by a spatially correlated noise.

\begin{definition}
A Baire space $\mathcal{B}$ is a separated topological space in which all
countable intersections of dense open subsets are dense also. A subset $%
\mathcal{A}$ of $\mathcal{B}$ is said to be meager (or a first category set
in the Baire sense), if it is contained in a countable union of closed
nowhere dense subsets of $\mathcal{B}$. The complement of a meager set is
called a comeager (or residual or a second category set).
\end{definition}

Let us introduce some notations.

\noindent For any $p\geq 2$, let $\mathcal{M}^{p}=\{u:\Omega \times \mathbb{R}%
_{+}\times \mathbb{R}^{d}\longrightarrow \mathbb{R}$, jointly continuous in time
and space such that for all $T>0$ and a compact subset $K$ of $\mathbb{R}^{d}$, $%
E\left\| u\right\| _{\gamma ,T,K}^{p}<+\infty \}$. Define a metric on $%
\mathcal{M}^{p}$ by: 
\[
d\left( u_{1},u_{2}\right) =\left( E\left\| u_{1}-u_{2}\right\| _{\gamma
,T,K}^{p}\right) ^{\frac{1}{p}}. 
\]
By using the Borel--Cantelli lemma, it is easy to see that $\left( \mathcal{M%
}^{p},d\right) $ is a complete metric space.\newline

\noindent Let $\mathcal{C}_{1}$ be the set of functions $b:\mathbb{R}_{+}\times 
\mathbb{R}^{d}\times \mathbb{R}\longrightarrow \mathbb{R}$ which are continuous and
bounded by $\kappa $. Define the metric $\rho _{1}$ on $\ \mathcal{C}_{1}$
as follows: 
\[
\rho _{1}\left( b_{1},b_{2}\right) =\sum_{n=1}^{+\infty }\frac{1}{2^{n}}%
\cdot \frac{\rho _{1,n}(b_{1}-b_{2})}{1+\rho _{1,n}(b_{1}-b_{2})},
\]
where 
\[
\rho _{1,n}(h)=\sup_{t\in [0,T]}\sup_{x\in \mathbb{R}^{d}}\sup_{\left|
r\right| \leq n}\left| h(t,x,r)\right| .
\]
Notice that the metric $\rho _{1}$ is compatible with the topology of
uniform convergence on compact subsets of $\mathbb{R}_{+}\times \mathbb{R}%
^{d}\times \mathbb{R}.$\newline

\noindent Let $\mathcal{C}_{2}$ the set of continuous $\kappa $--bounded
functions $\sigma:\mathbb{R}_{+}\times \mathbb{R}^{d}\times \mathbb{R}\longrightarrow 
\mathbb{R}$ with the corresponding metric $\rho _{2}$. Clearly the space $%
\mathcal{R}=\mathcal{C}_{1}\times \mathcal{C}_{2} $ endowed with the product
metric $\rho $ is a complete metric space, hence $(\mathcal{R},\rho )$ is a
Baire space.\newline

\noindent Let $\mathcal{R}_{Lip}$ be the subset of $\mathcal{R}$ consisting
of functions which are continuous and satisfy \textbf{(L)} and \textbf{(LG)}.

\begin{proposition}
\label{pro11}The space $\mathcal{R}_{Lip}$ is a dense subset in $\mathcal{R}$%
.
\end{proposition}

\noindent\textbf{Proof.} By truncation and regularization arguments. \hfill $%
\Box$

\subsection{The existence is generic}

We denote by $\mathcal{R}_{e}$ the subset of functions $\sigma ,b$ in $%
\mathcal{R}$ for which equation (\ref{eq0}) has a, not necessarily unique,
solution and by $\mathcal{R}_{u}$ the subset of $\mathcal{R}$ which consists
to all functions $\sigma ,b$ for which equation (\ref{eq0}) has a unique
solution.

\begin{theorem}
\label{th11}$\mathcal{R}_{u}$ is a residual set in the Baire space $(%
\mathcal{R},\rho )$.
\end{theorem}

\noindent To prove this Theorem we need some lemmas.

\begin{lemma}
\label{lem22}Let $\sigma ,b$ be elements of $\mathcal{R}_{Lip}$. Let $%
(\sigma _{n},b_{n})_{n\in \mathbb{N}}$ be a sequence in $\mathcal{R}_{e}$. We
assume that 
\[
\rho \left[ (\sigma _{n},b_{n}),(\sigma ,b)\right] \longrightarrow 0\quad %
\hbox {as}\quad n\longrightarrow \infty .
\]
Then, under $\mathbf{(R.1)}$ and $\mathbf{(R.2)}$, $u^{\sigma _{n},b_{n}}$
converges to $u^{\sigma ,b}$ in $(\mathcal{M}^{p},d)$ as $n\rightarrow
+\infty $.
\end{lemma}

\noindent \textbf{Proof.} To lighten the notation we write $u$ for $%
u^{\sigma ,b}$ and $u_{n}$ for $u^{\sigma _{n},b_{n}}$. It is clear that $%
u_{n}-u$ satisfies the property P$_{1}$ of Lemma \ref{lem1}. \newline

\noindent Let us put $\varphi _{n}(t,x)=E|u_{n}(t,x)-u(t,x)|^{p}$ and $\phi
_{n}(t)=\sup_{x\in \mathbb{R}^{d}}\varphi _{n}(t,x)$. Hence 
\begin{eqnarray}
\varphi _{n}(t,x) &\leq &c_{p}(\nu _{t})^{\frac{p}{2}-1}\left[
\int_{0}^{t}ds\left( \sup_{y\in \mathbb{R}^{d}}E|\sigma _{n}(u_{n})(s,y)-\sigma
(u)(s,y)|^{p}\right) J(t-s),\right.   \nonumber \\
&&+\left. \int_{0}^{t}ds\left( \sup_{y\in \mathbb{R}%
^{d}}E|b_{n}(u_{n})(s,y)-b(u)(s,y)|^{p}\right) J(t-s)\right] ,  \nonumber
\end{eqnarray}
where $J$ and $\mu $ are defined by (\ref{Jmu}). Fix $T>0$ and define for
each integer $N$ 
\[
\tau _{n}^{N}=\inf \left\{ t\geq 0:\sup_{x\in \mathbb{R}^{d}}|u_{n}(t,x)|\geq
N\right\} .
\]
The uniform boundedness of $\sigma _{n}$, $b_{n}$, $\sigma $ and $b$ imply
that 
\[
\lim_{N\rightarrow +\infty }P(\left\| u_{n}\right\| _{T,\infty }\geq N)=0,
\]
uniformly with respect to $n$, whence 
\[
\lim_{N\rightarrow +\infty }P(\tau _{n}^{N}>T)=1,
\]
uniformly with respect to $n$.\newline

Let $\varepsilon >0$ be fixed and choose $N(\varepsilon )$ a natural number
such that $N(\varepsilon )>T$ and $P(\tau _{n}^{N(\varepsilon
)}>T)>1-\varepsilon $, for all $n$. We have 
\begin{eqnarray*}
\lefteqn{\sup_{y\in \mathbb{R}^{d}}E|\sigma _{n}(u_{n})(s,y)-\sigma
(u)(s,y)|^{p}} \\
&\leq &c_{p}\sup_{y\in \mathbb{R}^{d}}E\left| \sigma _{n}(u_{n})(s,y)-\sigma
(u_{n})(s,y)\right| ^{p} \\
&&+c_{p}\sup_{y\in \mathbb{R}^{d}}E\left| \sigma (u_{n})(s,y)-\sigma
(u)(s,y)\right| ^{p} \\
&\leq &c_{p}\left( \sup_{y\in \mathbb{R}^{d}}E\left[ \left| \sigma
_{n}(u_{n})(s,y)-\sigma (u_{n})(s,y)\right| ^{p}\mathbf{1}_{\{s\leq \tau
_{n}^{N(\varepsilon )}\}}\right] \right.  \\
&&+\sup_{y\in \mathbb{R}^{d}}E\left[ \left| \sigma _{n}(u_{n})(s,y)-\sigma
(u_{n})(s,y)\right| ^{p}\mathbf{1}_{\{s>\tau _{n}^{N(\varepsilon )}\}}\right]
\\
&&\left. +\sup_{y\in \mathbb{R}^{d}}E\left| \sigma (u_{n})(s,y)-\sigma
(u)(s,y)\right| ^{p}\right)  \\
&\leq &c_{p}\sup_{t\in [0,T]}\sup_{x\in \mathbb{R}^{d}}\sup_{|r|\leq
N(\varepsilon )}|\sigma _{n}(t,x,r)-\sigma (t,x,r)|^{p} \\
&&+c_{p}P\left( \tau _{n}^{N(\varepsilon )}<s\right) +c_{p}\phi _{n}(s) \\
&\leq &c_{p}\left( \rho _{2,N(\varepsilon )}^{p}(\sigma _{n}-\sigma
)+P\left( \tau _{n}^{N(\varepsilon )}<s\right) +\phi _{n}(s)\right) ,
\end{eqnarray*}
where we have used the Lipschitz property of $\sigma $.\newline

\noindent With similar arguments we obtain 
\begin{eqnarray*}
\lefteqn{\sup_{y\in \mathbb{R}^{d}}E|b_{n}(u_{n})(s,y)-b(u)(s,y)|^{p}\,\,\,\,\,%
\,} \\
&\leq &c_{p}\left( \rho _{1,N(\varepsilon )}^{p}(b_{n}-b)+P(\tau
_{n}^{N(\varepsilon )}<s)+\phi _{n}(s)\right) .
\end{eqnarray*}

Set $\nu _{T}^{\ast }=\sup_{t\leq T}\nu _{t}$. Then for $0\leq t\leq T$, we
have 
\begin{eqnarray*}
\phi _{n}(t) &\leq &c_{p}\,\nu _{t}^{\frac{p}{2}-1}\left[ \,\int_{0}^{t}%
\left( \rho _{2,N(\varepsilon )}^{p}(\sigma _{n}-\sigma )+\rho
_{1,N(\varepsilon )}^{p}(b_{n}-b)\right. \right. \\
&&\quad \quad \qquad +\left. \left. P\left( \tau _{n}^{N(\varepsilon
)}<s\right) +\phi _{n}(s)\right) J(t-s)ds\right] \\
&\leq &c_{p}\,\nu _{t}^{\frac{p}{2}-1}\left[ \nu _{t}\left( \rho
_{2,N(\varepsilon )}^{p}(\sigma _{n}-\sigma )+\rho _{1,N(\varepsilon
)}^{p}(b_{n}-b)\right) \right. \\
&&\quad \quad \qquad \left. +\,\nu _{t}P\left( \tau _{n}^{N(\varepsilon
)}<t\right) +\int_{0}^{t}\phi _{n}(s)J(t-s)ds\right] \\
&\leq &c_{p}\,(\nu _{T}^{\ast })^{\frac{p}{2}-1}\left[ \nu _{T}^{\ast }(\rho
_{2,N(\varepsilon )}^{p}(\sigma _{n}-\sigma )+\rho _{1,N(\varepsilon
)}^{p}(b_{n}-b))\right. \\
&&\quad \quad \qquad \left. +\,\nu _{T}^{\ast }P\left( \tau
_{n}^{N(\varepsilon )}<t\right) +\int_{0}^{t}\phi _{n}(s)J(t-s)ds\right] .
\end{eqnarray*}
Hence, using Lemma \ref{lemGG} we get 
\[
\sup_{0\leq t\leq T}\phi _{n}(t)\leq c_{p,T}\left( \rho _{2,N(\varepsilon
)}^{p}(\sigma _{n}-\sigma )+\rho _{1,N(\varepsilon )}^{p}(b_{n}-b)+P(\tau
_{n}^{N(\varepsilon )}<T)\right) . 
\]
In view of the convergence of $\sigma _{n}$, $b_{n}$ to $\sigma $, $b$ it
follows that 
\[
\lim_{n\rightarrow +\infty }\left[ \rho _{2,N(\varepsilon )}^{p}(\sigma
_{n}-\sigma )+\rho _{1,N(\varepsilon )}^{p}(b_{n}-b)\right] =0. 
\]
Moreover, using the definition of $N(\varepsilon )$, and the above inequality, it
follows that 
\[
\lim_{n\rightarrow +\infty }\sup_{t\leq T}\phi _{n}(t)\leq
c_{p,T}\varepsilon , 
\]
and therefore, since $\varepsilon $ and $T$ are arbitrary, Lemma \ref{lem22}
follows now by tightness criterion which is assured by \textbf{(R.1}) and 
\textbf{(R.2)}.

Now, we define the oscillation function $\theta:\mathcal{R}\longrightarrow 
\mathbb{R}_{+}$ as 
\begin{eqnarray*}
\theta (\sigma ,b) &=&\lim_{\delta \rightarrow 0}\left[ \sup \left\{
d(u^{\sigma _{1},b_{1}},u^{\sigma _{2},b_{2}})\,\,\mbox{{s.t.}}\,\,(\sigma
_{i},b_{i})\in \mathcal{R}_{Lip}\mbox{, }i=1,2\right. \right. \\
&&\;\;\;\;\mbox{\ and }\left. \left. \max \left\{ \rho \left[ (\sigma
,b);(\sigma _{1},b_{1})\right] ,\rho \left[ (\sigma ,b);(\sigma _{2},b_{2})%
\right] \right\} <\delta \right\} \right]
\end{eqnarray*}
Then, we have the following:

\begin{lemma}
\label{lem33}Assume that $\mathbf{(R.1})$ and $\mathbf{(R.2)}$ hold$\mathbf{.%
}$\newline
(i) If $(\sigma ,b)$ belongs to $\mathcal{R}_{Lip}$ then $\theta (\sigma
,b)=0$.\newline
(ii) The function $\theta $ is upper semicontinuous on $\mathcal{R}$.\newline
(iii) If $\theta (\sigma ,b)=0$ for $\sigma ,b$ in $\mathcal{R}$, then
equation $(\ref{eq1})$ has at least one solution in $\mathcal{M}^{p}$.
\end{lemma}

\begin{remark}
The assertion $(iii)$ of Lemma \ref{lem33} is a sufficient condition to
ensure existence of solutions of equation $(\ref{eq1})$.
\end{remark}

\noindent \textbf{Proof of Lemma \ref{lem33}.} The assertion $(i)$ is no
more than an immediate consequence of Lemma \ref{lem22}.\newline
Proof of $(ii)$. By the definition of $\theta$ we can write: 
\begin{eqnarray*}
\theta (\sigma ,b) &=&\lim_{\delta \rightarrow 0}\theta_{\delta} (\sigma ,b)
\end{eqnarray*}
where 
\begin{eqnarray*}
\theta_{\delta} (\sigma ,b)&=&\left[ \sup \left\{ d(u^{\sigma
_{1},b_{1}},u^{\sigma _{2},b_{2}})\,\,\mbox{{s.t.}}\,\,(\sigma
_{i},b_{i})\in \mathcal{R}_{Lip}\mbox{, }i=1,2\right. \right. \\
&&\;\;\;\;\mbox{\ and }\left. \left. \max \left\{ \rho \left[ (\sigma
,b);(\sigma _{1},b_{1})\right] ,\rho \left[ (\sigma ,b);(\sigma _{2},b_{2})%
\right] \right\} <\delta \right\} \right].
\end{eqnarray*}
Let $\eta >0$ and $(\sigma ,b)$ in $\mathcal{R}$. It is not difficult to see
that for every $\delta > 0$, $\theta_{\delta} (\sigma ,b)\leq
\theta_{\delta+\eta} (\sigma_1 ,b_1)$ for each $(\sigma_1 ,b_1)$ in $%
\mathcal{R}$ such that $\rho[(\sigma ,b);(\sigma_1 ,b_1)] < \eta.$ Which
implies that the mapping $\theta$ is upper semicontinuous.

\noindent We now turn to the proof of $(iii)$. Let $(\sigma ,b)\in \mathcal{R%
}$. Since $\theta (\sigma ,b)=0$, then there exists a decreasing sequence of
strictly positive numbers $\delta _{n}$ $(\delta _{n}\searrow 0)$ such that 
\begin{eqnarray}
&&\sup \left\{ d(u^{1},u^{2})\,\,\mbox{such
that}\,\,(\sigma _{i},b_{i}))\in \mathcal{R}_{Lip},\;i=1,2\,\,\right.
\label{enq3} \\
&&\;\;\;\left. \mbox{ and }\max \left\{ \rho \left[ (\sigma ,b);(\sigma
_{1},b_{1})\right] ;\,\rho \left[ (\sigma ,b);(\sigma _{2},b_{2})\right]
\right\} <\delta \right\} <\frac{1}{n}  \nonumber
\end{eqnarray}
But then, by Proposition \ref{pro11}, for each $n\in \mathbb{N}^{\ast }$, there
exists a $(\sigma ^{n},b^{n})\in \mathcal{R}_{Lip}$ such that $\rho [(\sigma
^{n},b^{n});(\sigma ,b)]<\delta _{n}$. Since $\delta _{n}$ decreases, it
follows from (\ref{enq3}) that $d(u_{n},u_{m})<\max (\frac{1}{m},\frac{1}{n})
$. Hence, $(u_{n})_{n\in \mathbb{N}}$ is a Cauchy sequence in the Banach space $%
(\mathcal{M}^{p},d)$. Thus, there exists $u\in \mathcal{M}^{p}$ such that 
\begin{equation}
\lim_{n\rightarrow \infty }d(u_{n},u)=0.  \label{enq4}
\end{equation}
Let us now check that $u$ satisfies equation (\ref{eq1}). Indeed, from (\ref
{enq4}), there exists a subsequence $n_{k}$ such that 
\begin{equation}
u_{n_{k}}(t,x)\quad \hbox{converges to}\quad u(t,x)\;\;\;P\mbox{--a.s. }\,\,%
\mbox{as}\,\,k\rightarrow +\infty .  \label{enq7}
\end{equation}
It remains now to prove that for each $t\in \mathbb{R}_{+}$, $I_{n}(t,x)$ and $%
J_{n}(t,x)$ converge in Probability to $0$ as $n$ goes to infinity, where 
\[
I_{n}(t,x)=\int_{0}^{t}ds\int_{\mathbb{R}^{d}}dyS(t-s,x,y)\left[
b^{n}(u_{n})(s,y)-b(u)(s,y)\right] 
\]
and 
\[
J_{n}(t,x)=\int_{0}^{t}\int_{\mathbb{R}^{d}}S(t-s,x,y)\left[ \sigma
^{n}(u_{n})(s,y)-\sigma (u)(s,y)\right] M(ds,dy). 
\]
We have 
\begin{eqnarray*}
E\left| J_{n}(t,x)\right| ^{p} &\leq &c_{p}\nu _{t}^{\frac{p}{2}%
-1}\int_{0}^{t}ds\left( \sup_{y\in \mathbb{R}^{d}}E|\sigma
^{n}(u_{n})(s,y)-\sigma (u)(s,y)|^{2}\right) J(t-s) \\
&\leq &\frac{c_{p}\nu _{t}^{\frac{p}{2}}}{N^{2}}+c_{p}\rho _{2,N}^{p}(\sigma
^{n}-\sigma )\nu _{t}^{\frac{p}{2}} \\
&&+c_{p}\nu _{t}^{\frac{p}{2}-1}\int_{0}^{t}ds\left( \sup_{y\in \mathbb{R}%
^{d}}E|\sigma (u_{n})(s,y)-\sigma (u)(s,y)|^{p}\right) J(t-s).
\end{eqnarray*}
Lemma \ref{lem22} implies that 
\[
\lim_{n,N\rightarrow +\infty }\frac{c_{p}}{N^{2}}\nu _{T}^{\frac{p}{2}%
}+c_{p}\rho _{2,N}^{p}(\sigma ^{n}-\sigma )\nu _{T}^{\frac{p}{2}}=0. 
\]
On the other hand, since $\sigma \in \mathcal{R}_{Lip}$ then (\ref{enq7})
implies that $\sigma (u_{n})$ converges to $\sigma (u)$ $dP\times dt\times
dx $--a.e. Hence, the Lebesgue dominated convergence theorem allows us to
deduce that 
\[
\lim_{n\rightarrow \infty }\int_{0}^{t}ds\left( \sup_{y\in \mathbb{R}%
^{d}}E|\sigma (u_{n})(s,y)-\sigma (u)(s,y)|^{p}\right) J(t-s)=0. 
\]
The proof of $I_{n}(t,x)$ can be carried out by similar argument. This
proves assertion (iii).\newline

\noindent \textbf{Proof of Theorem \ref{th11}.} Lemma \ref{lem22} and
assertions (i) and (ii) of Lemma \ref{lem33} imply that for each natural
number $n$, the set $\mathcal{H}_{n}=\{(\sigma ,b)\in \mathcal{R}:\theta
(\sigma ,b)<\frac{1}{n}\}$ is a dense open subset of $(\mathcal{R},\rho )$.
Then by the Baire categories theorem the set $\mathcal{H}=\cap _{n\in \mathbb{N}%
^{\ast }}\mathcal{H}_{n}$ is a dense $G_{\delta }$ subset of the Baire space 
$(\mathcal{R},\rho )$. Moreover, if $(\sigma ,b)\in \mathcal{H}$ then Lemma 
\ref{lem33} (iii) implies that the corresponding equation (\ref{eq1}) has
one solution. Hence $\mathcal{H}\subset \mathcal{R}_{e}$. This implies that $%
\mathcal{R}_{e}$ is a residual subset in $(\mathcal{R},\rho )$.\newline
To prove that $\mathcal{R}_{u}$ is residual, we define the function $%
\widetilde{\theta }:\mathcal{H}\longrightarrow \mathbb{R}_{+}$ as follows, 
\[
\widetilde{\theta }(\sigma ,b)=\sup \left\{ d(u_{1}^{\sigma
,b},u_{2}^{\sigma ,b}):u_{i}^{\sigma ,b}\hbox { is a solution to equation }(%
\ref{eq1}),\;i=1,2\right\} 
\]
and for each $n\in \mathbb{N}^{\ast }$ we put $\overline{\mathcal{H}}%
_{n}=\{(\sigma ,b)\in \mathcal{H}:\widetilde{\theta }(\sigma ,b)<\frac{1}{n}%
\}$. By using Lemma \ref{lem22} we see, as in the proof of Lemma \ref{lem33}
(ii), that the function $\widetilde{\theta }$ is upper semicontinuous on $%
\mathcal{R}$. This implies that each $\overline{\mathcal{G}}_{n}$
contains the intersection of $\mathcal{H}$ and a dense open subset of $(%
\mathcal{R},\rho )$. Thus the set $\overline{\mathcal{H}}=\cap _{n\in \mathbb{N}%
^{\ast }}\overline{\mathcal{H}}_{n}$ contains a dense $G_{\delta }$ subset
of the Baire space $(\mathcal{R},\rho )$. Hence it is residual in $(\mathcal{%
R},\rho )$. Finally, if $(\sigma ,b)\in \overline{\mathcal{H}}$ then the
corresponding equation (\ref{eq1}) has a unique solution. Thus $\overline{%
\mathcal{H}}\subset \mathcal{R}_{u}$. Theorem \ref{th11} follows.

\subsection{Continuous dependence on the coefficients}

For a given $\sigma ,b\in \mathcal{R}$ we denote by $\Phi (\sigma
,b)=u^{\sigma ,b}$ the solution of Eq($\sigma ,b$) when it exists.

\begin{theorem}
\label{th22} There exists a second category set $\mathcal{R}_{2}$ such that
the map $\Phi:\mathcal{R}_{2}\longrightarrow \mathcal{M}^{p}$ given by $\Phi
(\sigma ,b)=u^{\sigma ,b}$ is well defined and continuous at each point of $%
\mathcal{R}_{2}$.
\end{theorem}

\noindent \textbf{Proof.} We shall show that $\Phi $ is continuous on $%
\overline{\mathcal{H}}$ (the dense $G_{\delta }$ set which has been defined
in the proof of Theorem \ref{th11}). Suppose the contrary. Then there exist $%
\sigma \in \overline{\mathcal{H}}$, $\varepsilon >0$ and a sequence $(\sigma
^{p})_{p}\subset \overline{\mathcal{H}}$ such that, 
\begin{equation}
\lim_{p\rightarrow +\infty }\rho \left[ (\sigma ^{p},b^{p});(\sigma ,b)%
\right] =0\;\;\mbox{and}\;\;d\left[ \Phi (\sigma ^{p},b^{p});\Phi (\sigma ,b)%
\right] \geq \varepsilon \mbox{ \ for each }p.  \label{enq9}
\end{equation}
Fix $n\in \mathbb{N}$ such that $\varepsilon <\frac{1}{n}$. Since $\overline{%
\mathcal{H}}\subset \mathcal{H}$ then there exists a decreasing sequence of
strictly positive numbers $\delta _{n}$ $(\delta _{n}\searrow 0)$ and a
sequence of functions $\sigma ^{n},b^{n}\in \mathcal{R}_{Lip}$ such that, 
\begin{equation}
\rho \left[ (\sigma ^{n},b^{n});(\sigma ,b)\right] <\delta _{n}\;\;\mbox{and}%
\;\;d\left[ \Phi (\sigma ^{n},b^{n}),\Phi (\sigma ,b)\right] <\frac{1}{n}.
\label{enq10}
\end{equation}
We choose $p$ large enough to have $\rho ((\sigma ^{p},b^{p});(\sigma
,b)]<\delta _{n}-\rho [(\sigma ^{n},b^{n});(\sigma ,b)]$ then we use (%
\ref{enq10}) to obtain $\rho [(\sigma ^{p},b^{p});(\sigma
^{n},b^{n})]<\delta _{n}.$ Hence by Lemma \ref{lem22} we have $d[[\Phi
(\sigma ^{p},b^{p});\Phi (\sigma ^{n},b^{n})]<\frac{1}{n}$. Thus 
\begin{eqnarray*}
d[\Phi (\sigma ^{p},b^{p});\Phi (\sigma ,b)] &\leq &d[\Phi (\sigma
^{p},b^{p});\Phi (\sigma ^{n},b^{n})] \\
&&+d[\Phi (\sigma ^{n},b^{n}),\Phi (\sigma ,b)] \\
&<&\frac{2}{n}
\end{eqnarray*}
which contradicts (\ref{enq9}). Theorem \ref{th22} is then proved.

\subsection{The Uniqueness is generic}

The main result of this section is the following

\begin{theorem}
\label{thgpu}The subset $\mathcal{R}_{pu}$ of $\mathcal{R}$ consisting of
those $\left( \sigma ,b\right) $ for which the property $\mathbf{(PU)}$
holds for $\mathit{Eq}(\sigma ,b)$ is a residual set.
\end{theorem}

\begin{lemma}
For each $\left( \sigma _{1},b_{1}\right) \in \mathcal{R}_{Lip}$ and $%
\varepsilon >0$ there exists $\delta \left( \varepsilon \right)
>0 $\ such that $\forall \;\left( \sigma ,b\right) \in B\left( \left( \sigma
_{1},b_{1}\right) ,\delta \right) $ and for every pair of solutions $u,v$
of $Eq(\sigma ,b)$ (defined on the same probability space, with
the same martingale measure), we have $d\left( u,v\right) <\varepsilon .$
\end{lemma}

\noindent \textbf{Proof}. Let $w$ be the unique strong solution of the SPDEs
Eq($\sigma _{1},b_{1}$) defined on the same probability space, and with
respect to the same martingale measure $M$. We have 
\[
d\left( u,v\right) \leq d\left( u,w\right) +d\left( w,v\right) , 
\]
the result follows from the continuity of $w$ with respect to the
coefficients.

\noindent \textbf{Proof of Theorem \ref{thgpu}.} We put $\mathcal{K}=\cap
_{k\geq 1}\cup _{\left( \sigma ,b\right) \in \mathcal{R}_{Lip}}B((\sigma
,b),\delta (\frac{1}{k}))$, the subset $\mathcal{K}$ is a $G_{\delta }$
dense subset in the Baire space $(\mathcal{R},\rho )$, and for every $\left(
\sigma ,b\right) \in \mathcal{K}$, the pathwise uniqueness holds for the
SPDEs Eq($\sigma ,b$). It follows that $\mathcal{R}_{pu}$ is a residual
subset in $\mathcal{R}$.

\begin{remark}
Gy\"{o}ngy (2001) has shown that $\mathcal{R}\setminus \mathcal{K}$ is not
empty.
\end{remark}

\section{Applications}

This section is devoted to study two examples for which Theorems \ref{th1}, \ref
{thip}, \ref{thips}, \ref{th11}, \ref{th22} and \ref{thgpu} can be applied.

\subsection{Stochastic heat equation}

Consider the following SPDEs 
\[
\left\{ 
\begin{array}{l}
\dfrac{\partial u}{\partial t}(t,x)=\dfrac{1}{2}\Delta u(t,x)+\sigma (u(t,x))%
{\dot{M}}(t,x)+b(u(t,x)) \\ 
\\ 
\;\;u(0,x)=0,\; t\in [0,T], \; x\in \mathbb{R}^{d}, \; d\geq 1.
\end{array}
\right. 
\]
Recall that we consider null initial conditions for the sake of simplicity.
The solution to this equation is given by 
\begin{eqnarray}
u(t,x) &=&\int_{0}^{t}\int_{\mathbb{R}^{d}}S_{h}^{d}(t-s,x-y)\sigma
(u(s,y))M(ds,dy)  \label{heat1} \\
&&  \nonumber \\
&&+\int_{0}^{t}ds\int_{\mathbb{R}^{d}}dyS_{h}^{d}(t-s,x-y)b(u(s,y)),  \nonumber
\end{eqnarray}
where $S_{h}^{d}(t,x)=(2\pi t)^{-d/2}\exp (\frac{-|x|^{2}}{2t})$ is the
fundamental solution to heat equation, with $d$--dimensional spatial
parameter.

\subsection{Stochastic wave equation}

This section deals with the SPDEs 
\[
\left\{ 
\begin{array}{l}
\dfrac{\partial ^{2}v}{\partial t^{2}}(t,x)=\Delta v(t,x)+b(v)(t,x)+\sigma
(v)(t,x){\dot{M}}(t,x) \\ 
\;\;v(0,x)=v_{0}(x)\;\mbox{ and}\; \dfrac{\partial v}{\partial t}(0,x)=0,\;
t\geq 0, x\in \mathbb{R}^{d} \; \mbox{ and} \; d\in \{1,2\}.
\end{array}
\right. 
\]
The solution to this equation is given by 
\begin{eqnarray}
v(t,x) &=&\int_{\mathbb{R}^{d}}W_{t}^{d}(x-y)v_{0}(y)dy+\int_{0}^{t}\int_{\mathbb{R%
}^{d}}W_{t-s}^{d}(x-y)\sigma (v)(s,y)M(ds,dy)  \nonumber \\
&&+\int_{0}^{t}\int_{\mathbb{R}^{d}}W_{t-s}^{d}(x-y)b(v)(s,y)dyds  \nonumber
\end{eqnarray}
where $W_{t}^{d}(x)$ is the fundamental solution to the wave equation with $%
d $--dimensional spatial parameter: 
\[
\left\{ 
\begin{array}{l}
W_{t}^{1}(x)=\dfrac{1}{2}\mathbf{1}_{\{|x|\leq t\}}\;\;\;\;\;\;\;\;\;\;\;%
\mbox{for \ }t\geq 0,\;x\in \mathbb{R} \\ 
\\ 
W_{t}^{2}(x)=\dfrac{1}{2\pi }\dfrac{1}{\sqrt{t^{2}-|x|^{2}}}\mathbf{1}%
_{\{|x|<t\}}\;\;\mbox{for \ }t\geq 0,\;x\in \mathbb{R}^{2}\mbox{.}
\end{array}
\right. 
\]

Let us now consider the stochastic wave equation in dimension $3$ whose
integral equation is given by 
\begin{eqnarray}
v(t,x) &=&\int_{\mathbb{R}^{3}}W_{t}^{3}(dy)v_{0}(x-y)+\int_{0}^{t}\int_{\mathbb{R}%
^{3}}W_{t-s}^{3}(x-y)\sigma (v(s,y))M(ds,dy)  \nonumber  \label{wave3} \\
&&+\int_{0}^{t}ds\int_{\mathbb{R}^{3}}b(v(t-s,x-y))W_{s}^{3}(dy),
\end{eqnarray}
where $W_{s}^{3}(dy)$ is the fundamental solution of the wave equation in $%
\mathbb{R}^{3}$. More precisely for each $s\in [0,T]$ $W_{s}^{3}=\frac{%
\Sigma _{s}}{4\pi s}$, where $\Sigma _{s}$ denotes the uniform surface
measure, with total mass $4\pi s^{2}$, on the sphere of radius $s$.\newline
It is proved by Dalang and Sanz--Sol\'{e} \cite{DaSa} that under the
following conditions:

\begin{enumerate}
\item  The covariance measure is absolutely continuous with respect to
Lebesgue measure with density given by $\Gamma (dx)=\varphi (x)|x|^{-\beta
}dx$, where $\varphi \in \mathcal{C}^{1}(\mathbb{R}^{3})$ bounded and positive
and $\triangledown \varphi \in \mathcal{C}_{b}^{\delta }(\mathbb{R}^{3})$ for
some $\beta \in ]0,2[$ and $\delta \in ]0,1[$.

\item  The support of the function $v_{0}$ is contained in the ball $%
B_{r_{0}}(0)$, for some $r_{0}>0.$

\item  The initial value function $v_{0}$ belongs to $\mathcal{C}^{2}(\mathbb{R}%
^{3})$ and $\Delta v_{0}$ is $\gamma $--H\"{o}lder continuous with $\gamma
\in ]0,1[$.
\end{enumerate}

The solution $v$ of the equation (\ref{wave3}) is $\alpha $--H\"{o}lder
continuous with $\alpha \in ]0,\gamma \wedge \frac{2-\beta }{2}\wedge \frac{%
1+\delta }{2}[$ jointly in $(t,x)$. \hfill $\Box$

\section{Appendix}

\begin{lemma}
\label{lemGG} ( \cite{CP} p. 226) Let $f$ be a positive function on $[0,T]$
such that 
\[
f(t)\leq h(t)+\int_{0}^{t}g(t-s)f(s)ds,\;\;t\in [0,T], 
\]
with $g\in L^{1}([0,T])$ and $h\in L^{p}([0,T])$, $p\geq 1$ both positive.
Then 
\[
f(t)\leq h(t)+\sum_{n=1}^{+\infty }(G^{n}h)(t)\mbox{ \ for all \
}t\in [0,T], 
\]
where, for $n\geq 1$, $G^{n}(h)=\int_0^t g_n(t-s)h(s)ds$ and $%
g_n(t)=\int_0^tg(t-s)g_{n-1}(s)ds$, $g_1(t)=g(t)$.\newline
In particular, if $h\equiv 0$ then $f\equiv 0.$
\end{lemma}

\begin{lemma}
\label{lem1}Let $\{Y_{n}(t,x)\}_{n\geq 1}$
be a sequence of processes indexed by $[0,T]\times \mathbb{R}^{d}$ such that  
\newline

\noindent $(\mathbf{P}_{1})$ For any $p\geq 2$ there exist $%
\,c_{p},\delta _{1},\delta _{2}>0$ such that for any $t_{1},t_{2}\in [0,T]$
and $x_{1},x_{2}\in \mathbb{R}^{d}$ 
\[
\sup_{n}E\left[ \left| Y_{n}(t_{2},x_{2})-Y_{n}(t_{1},x_{1})\right| ^{p}%
\right] \leq c_{p}\left( \left| t_{2}-t_{1}\right| ^{\delta _{1}p}+\left\|
x_{2}-x_{1}\right\| ^{\delta _{2}p}\right). 
\]

\noindent $(\mathbf{P}_{2})$ For every $(t,x)\in [0,T]\times \mathbb{R}%
^{d}$ and $p\geq 2$ 
\[
\lim_{n\rightarrow +\infty }E\left[ \left| Y_{n}(t,x)\right| ^{p}\right] =0. 
\]
Then for any $\gamma _{1}\in (0,\delta _{1})$ and $\gamma _{2}\in (0,\delta
_{2})$ 
\[
\lim_{n\rightarrow +\infty }E\left[ \left\| Y_{n}\right\| _{\gamma ,T,K}^{p}%
\right] =0, 
\]
where $K$ is a compact subset of $\mathbb{R}^{d}$.
\end{lemma}

\noindent \textbf{Proof.} The proof of this lemma can be found in Millet and
Sanz--Sol\'{e} \cite{MS}.\hfill $\Box $

\begin{lemma}
Let $f, h$ be two functions defined on $\mathbb{R}$ and $\mu$ a positive
measure such that $f \cdot h \in L^1(\mu)$. Then, for all $q>1$, we have: 
\begin{eqnarray}
\left\vert\int f\cdot \vert h \vert d\mu \right\vert^q\le \left(\int \vert
f\vert^q\cdot \vert h\vert d\mu\right) \left(\int \vert h\vert
d\mu\right)^{q-1}.  \nonumber
\end{eqnarray}
\end{lemma}

\noindent \textbf{Proof.} Set $\nu =|h|d\mu $, then the result follows from
the H\"{o}lder's inequality applied to $\int fd\nu $. \hfill $\Box $
\smallskip 

\noindent The following elementary Lemma is an extension of Gronwall's Lemma
akin to lemma 3.3 established in Walsh \cite{W}.

\begin{lemma}
\label{gronwall} Let $\theta > 0$.  Let $(f_n, n\in \mathbb{N})$ be a sequence
of non-negative functions on $[0,T]$ and $\alpha, \beta$ be non-negative
real numbers such that for $0\le t\le T,\, n\geq 1$ 
\begin{eqnarray}
f_n(t)\le \alpha+ \int_0^t \beta\,f_{n-1}(s) (t-s)^{\theta -1}ds.
\label{gronvol}
\end{eqnarray}
If $\sup_{0\le t \le T}f_0(t)=M, $ then for $n\geq 1$, 
\begin{eqnarray}
f_n(t)\le \frac{1}{2} \left(\alpha + \alpha \exp \left(\frac{2 \beta
t^{\theta}}{\theta}\right) + \frac{M}{n!} \left( \frac{2 \beta t^\theta}{%
\theta}\right)^n\right).  \nonumber
\end{eqnarray}
In particular, $\sup_{n \geq 0}\sup_{0\le t \le T}f_n(t)<\infty $, and if $%
\alpha=0$, then $\sum_{n \geq 0} f_n(t)$ converges uniformly on $[0,T]$.
\end{lemma}

\noindent \textbf{Proof.} Let us prove by induction that, for $n\geq 1$, 
\begin{eqnarray}\label{induc}
f_{n}(t)\leq \alpha \left( 1+\sum_{1\leq k\leq n-1}\frac{2^{k-1}}{k!}\left( 
\frac{\beta t^{\theta }}{\theta }\right) ^{k}\right) +M\cdot \frac{2^{n-1}}{%
n!}\left( \frac{\beta t^{\theta }}{\theta }\right) ^{n}.
\end{eqnarray}
The initial step is readily checked: 
\[
f_{1}(t)\leq \alpha +\int_{0}^{t}\beta \,M(t-s)^{\theta -1}ds=\alpha +M\,%
\frac{\beta t^{\theta }}{\theta }.
\]
Now, since (\ref{gronvol}) is true, we have 
\[
f_{n}(t)\leq \alpha +\int_{0}^{t}\beta \left( \alpha +\alpha \sum_{1\leq
k\leq n-2}\frac{2^{k-1}}{k!}\left( \frac{\beta s^{\theta }}{\theta }\right)
^{k}M\,\frac{2^{n-2}}{(n-1)!}\left( \frac{\beta s^{\theta }}{\theta }\right)
^{n-1}\right) (t-s)^{\theta -1}\,ds.
\]
Consider 
\[
\int_{0}^{t}s^{k\theta }(t-s)^{\theta -1}ds\leq
\int_{0}^{t/2}(t-s)^{(k+1)\theta -1}ds+\int_{t/2}^{t}s^{(k+1)\theta -1}ds.
\]
Hence we may bound 
\[
\int_{0}^{t}s^{k\theta }(t-s)^{\theta -1}ds\leq 2\frac{t^{(k+1)\theta }}{%
(k+1)\theta }.
\]
Summation over $k$ brings (\ref{induc}). \hfill $\Box $

\end{document}